\documentclass[journal]{IEEEtran}
\pdfoutput=1
\usepackage{amsfonts}
\usepackage{amsthm,amsmath, amssymb}
\usepackage{graphicx}
\usepackage{verbatim}
\usepackage{bbm}
\usepackage{rotating}

\usepackage{algorithmic}
\usepackage{algorithm}

\newcommand{\bo}[1]{\boldsymbol{#1}}

\newcommand{\cov}{\mbox{Cov}}
\newcommand{\E}{\mbox{E}}

\newtheorem{theorem}{Theorem}
\newtheorem{definition}{Definition}
\newtheorem{corollary}{Corollary}
\newtheorem{lemma}{Lemma}

\newcommand{\bmat}{\begin{pmatrix}}
\newcommand{\emat}{\end{pmatrix}}
\usepackage{color}

\graphicspath{ {./Figures/} }

\begin{document}
\title{The squared symmetric FastICA estimator}

\author{Jari~Miettinen, Klaus~Nordhausen, Hannu~Oja, Sara~Taskinen and Joni Virta% <-this % stops a space
\thanks{J. Miettinen and S. Taskinen are with the Department of Mathematics and Statistics,
University of Jyvaskyla, Jyvaskyla, FIN-40014, Finland
(e-mail: jari.p.miettinen@jyu.fi).}% <-this % stops a space
\thanks{K. Nordhausen, H. Oja and J. Virta are with the Department of Mathematics and Statistics, University of
Turku, Turku, FIN-20014, Finland.}% <-this % stops a space
}


\markboth{}
{Miettinen \MakeLowercase{\textit{et al.}}: The squared symmetric FastICA estimator}



\maketitle

\begin{abstract}
In this paper we study the theoretical properties of the deflation-based FastICA method, the original symmetric FastICA
method, and a modified symmetric FastICA method, here called the squared symmetric FastICA. This modification is obtained by replacing the absolute values in the FastICA objective function by their squares. In the deflation-based case this replacement has no effect on the estimate
 since the maximization problem stays the same. However, in the symmetric
case a novel estimate with unknown properties is obtained. In the paper we review the classic deflation-based
and symmetric FastICA approaches and contrast these with the new squared symmetric version of FastICA. We find the estimating equations and derive the asymptotical properties of the squared symmetric FastICA estimator with an arbitrary
choice of nonlinearity. Asymptotic variances of the unmixing matrix estimates are then used to compare their efficiencies
for large sample sizes showing that the squared symmetric FastICA estimator outperforms the
other two estimators in a wide variety of situations.
%{\color{red} The abstract must be a concise yet comprehensive reflection of what is in your article.
%The abstract must be self-contained, without abbreviations, footnotes, displayed equations, or references.
%The abstract must be between 150-250 words.}
\end{abstract}


% \fi
\IEEEpeerreviewmaketitle

\section{Introduction}

We assume that a $p$-variate random vector $\bo x=(x_1,\dots,x_p)^T$ follows the basic independent
component (IC) model, that is, the components of $\bo x$ are linear mixtures of $p$ mutually
independent latent variables in $\bo z=(z_1,\dots,z_p)^T$. The model can then be written as
\begin{align}
\label{ICmodel}
\bo x=\bo\mu+\bo\Omega\bo z,
\end{align}
where $\bo\mu$ is a location shift and $\bo\Omega$ is a full-rank $p\times p$ mixing matrix.
In independent component analysis (ICA), parameter $\bo\mu$ is usually regarded as a nuisance
parameter as the main interest is to find, using a random sample $\bo X=(\bo x_1,\dots,\bo x_n)$ from the
distribution of $\bo x$, an estimate for an unmixing matrix $\bo\Gamma$ such
that $\bo\Gamma\bo x$ has independent
components~\cite{HyvarinenKarhunenOja:2001,CichockiAmari:2002,ComonJutten:2010}.
Note that versions of \eqref{ICmodel} also exist where the dimension of $\bo z$ is larger than that of $\bo x$ (the \textit{underdetermined case}) or the other way around (the \textit{overdetermined case}), in the latter of which we can simply apply a dimension reduction method at first stage.
%{\color{blue} Notice that ICA can also be used in situations where the dimension of $\bo x$ is larger than that of
%$\bo z$ by applying a dimension reduction method at first stage.}
 In this paper we, however, restrict
to the case where $\bo x$ and $\bo z$ are of the same dimension.

The IC model~(\ref{ICmodel}) is a semiparametric model in the sense that the marginal distributions
of the components $z_1,\dots,z_p$ are unspecified. However, some assumptions on $\bo z$ are needed
in order to fix the model: For identifiability of $\bo\Omega$, we need to assume that
\begin{itemize}
\item[(A1)] at most one of the components $z_1,\dots,z_p$ is gaussian~\cite{Theis:2004}.
\end{itemize}
Nevertheless, $\bo\mu$, $\bo\Omega$ and $\bo z$ are still confounded and the mixing matrix
$\bo\Omega$ can be identified only up to the order, the signs, and heterogenous multiplication of its columns.
To fix $\bo\mu$ and  the scales of the columns of $\bo\Omega$ we further assume that
\begin{itemize}
\item[(A2)] $\E(z_i)=0$ and $\E(z_i^2)=1$ for $i=1,\dots,p$.
\end{itemize}
After these assumptions, the order and signs of the columns of $\bo\Omega$ still remain unidentified. 
For practical data analysis, this is, however, often sufficient.
The impact of the component order on asymptotics is further discussed in Section~\ref{asymptotics:sec}.

The solutions to the ICA problem are often formulated as algorithms with two steps. The first step is to whiten the
data, and the second step is to find an orthogonal matrix that rotates the whitened data to independent components.
In the following we formulate such an algorithm at the population level using the random variable $\bo x$:
Let $\bo S(F_{\bo x})=\bo{Cov}(\bo x)$ denote the covariance matrix of a random vector $\bo x$, 
where $F_{\bo x}$ denotes the cumulative distribution function $\bo x$, and write
$\bo x_{st}=\bo S^{-1/2}(F_{\bo x})(\bo x-\E(\bo x))$ for the standardized (whitened) random vector.
Here the square root matrix $\bo S^{-1/2}$ is chosen to be symmetric.
The aim of the second step is to find the rows of an orthogonal matrix
$\bo U=(\bo u_1,\dots,\bo u_p)^T$, either one by one (deflation-based approach) or
simultaneously (symmetric approach). The symmetric version of the famous FastICA
algorithm~\cite{Hyvarinen:1999} finds the orthogonal matrix $\bo U$, which maximizes
a measure of non-Gaussianity for the rotated components,
$$
\sum_{j=1}^p|\E[G(\bo u_j^T\bo x_{st})]|,
$$
where $G$ is a twice continuously differentiable, nonlinear and nonquadratic function (see
Section~\ref{Gfunction:sec} for more details). 

In this paper we replace the absolute values by their squares and consider the objective function
$$
\sum_{j=1}^p \left(\E[G(\bo u_j^T\bo x_{st})]\right)^2,
$$
as suggested
in~\cite{VirtaNordhausenOja:2015}, where the squared symmetric FastICA estimates based on convex combinations of the third and fourth squared cumulants were
studied in detail. Notice that replacing the absolute values by their squares in the objective functions has been  mentioned
in~\cite{Hyvarinen:1999} and~\cite[Section 6]{ComonJutten:2010}, but the idea was never carried  further.
In Section~\ref{functionals:sec} we formulate unmixing matrix functionals based on the two
symmetric approaches and the deflation-based approach.
Some statistical properties of the old estimators are recalled in Section~\ref{asymptotics:sec},
and the corresponding results of squared symmetric FastICA are derived for the first time for general
function $G$. The efficiencies of the three estimators are compared in
Section~\ref{comparisons:sec} using both asymptotic results and simulations.

\section{FastICA functionals}
\label{functionals:sec}

In this section we give formal definitions of three different, two old and one new, FastICA unmixing matrix functionals
with corresponding estimating equations and algorithms for their computation. The formal definition of the squared symmetric FastICA functional is new.
The conditions for function $G$ that ensure the consistency of the estimates is also discussed.
%the corresponding estimating equations and the algorithms. The estimating equations of the symmetric
%and the squared symmetric FastICA as well as the algorithm for the squared symmetric FastICA are given
%here for the first time. The section ends with the discussion on the conditions for function $G$ to
%ensure the consistency of the estimators.

\subsection{IC functionals}

Let again $F_{\bo x}$ denote the cumulative distribution function of a random vector $\bo x$ obeying
the IC model~(\ref{ICmodel}), and write $\bo\Gamma(F_{\bo x})$ for the value of an unmixing matrix
functional at the distribution $F_{\bo x}$. Due to the ambiguity in model~(\ref{ICmodel}), it is
natural to require that the separation result $\bo\Gamma(F_{\bo x})\bo x=\bo\Gamma(F_{\bo z})\bo z $ does not depend on
$\bo\mu$ and $\bo\Omega$ and the choice of $\bo z$ in the model specification. This is formalized in the following.
\begin{definition}\label{ICfunc}
The $p\times p$ matrix-valued functional $\bo\Gamma(F_{\bo x})$ is said to be
independent component (IC) functional  if
\begin{enumerate}
\item  $\bo\Gamma(F_{\bo x})\bo x$ has independent components for all $\bo x$ in the IC model (\ref{ICmodel}), and
\item  $\bo\Gamma(F_{\bo x})$ is affine-equivariant in the sense that
$\bo\Gamma(F_{\bo A\bo x+\bo b})=\bo\Gamma(F_{\bo x})\bo A^{-1}$ for all nonsingular $p \times p$
full-rank matrices $\bo A$, for all $p$-vectors $\bo b$ and for all $\bo x$ (even beyond the IC model).
\end{enumerate}
\end{definition}
The condition $\bo\Gamma(F_{\bo A\bo x+\bo b})=\bo\Gamma(F_{\bo x})\bo A^{-1}$
can be relaxed to be true only up to permutations and sign changes of their rows.
The corresponding sample version $\hat{\bo\Gamma}=\bo\Gamma(\bo X)$ is obtained when the IC functional
is applied to the empirical distribution function of $\bo X=(\bo x_1,\dots,\bo x_n)$. Naturally,
the estimator is then also affine equivariant in the sense that
$\bo\Gamma(\bo A\bo X+\bo b\bo 1_n^T)=\bo\Gamma(\bo X)\bo A^{-1}$ for all nonsingular $p \times p$ full-rank
matrices $\bo A$ and for all $p$-vectors $\bo b$.

The rest of this section focuses on three specific FastICA functionals. For recent overviews of FastICA and its
variants see also \cite{KoldovskyTichavsky:2015} and \cite{Wei:2015}.

\subsection{Deflation-based approach}

Deflation-based FastICA functional is based on the algorithm proposed in~\cite{HyvarinenOja:1997}
and~\cite{Hyvarinen:1999}. In deflation-based FastICA method the rows of an unmixing matrix are
extracted one after another. The method can thus be used in situations where only the few
most important components are needed. The statistical properties of the deflation-based method were
studied in~\cite{Ollila:2009} and \cite{Ollila:2010}, where the influence functions and limiting variances and
covariances of the rows of unmixing matrix were derived.

Assume now that $\bo x$ is an observation from an IC model~(\ref{ICmodel}) with mean vector
$\bo\mu=\E(\bo x)$ and covariance matrix $\bo S=\bo{Cov}(\bo x)$. In deflation-based FastICA,
the unmixing matrix $\bo\Gamma=(\bo\gamma_1,\dots,\bo\gamma_p)^T$ is estimated so that after finding
$\bo\gamma_1,\dots,\bo\gamma_{j-1}$, the $j$th row vector $\bo\gamma_j$
maximizes a measure of non-Gaussianity
 \[
|\E[G(\bo\gamma_j^T(\bo x-\E(\bo x)))]|
\]
under the constraints $\bo\gamma^T_l\bo S\bo\gamma_j=\delta_{lj}$, $l=1,\dots,j$, where
$\delta_{lj}$ is the Kronecker  delta $\delta_{lj}=1$ $(0)$ as $l=j$ $(l\neq j)$.
The requirements for the function $G$ and the conventional choices of it are discussed in Section~\ref{Gfunction:sec}.

%The $j$th FastICA functional $\bo\gamma_j=\bo\gamma_j(F_{\bo x})$ optimizes the Lagrangian function
%\begin{equation*}
%\begin{split}
%L(\bo\gamma_j,\bo\theta_j)=&|\E[G(\bo\gamma_j^T(\bo x-\E(\bo x)))]|-\frac{\theta_{jj}}{2}(\bo\gamma^T_j\bo S\bo\gamma_j-1) \\
% &-\sum_{i=1}^{j-1}\theta_{ji}\bo\gamma^T_i\bo S\bo\gamma_j,
%\end{split}
%\end{equation*}
%where $\bo\theta_j=(\theta_{j1},\dots,\theta_{jj})^T$ are the Lagrangian multipliers. By differentiating
%$L(\bo\gamma_j,\bo\theta_j)$ with respect to $\bo\gamma_j$, setting the derivative vector to zero and
%solving the Lagrangian multipliers, it is then easy to show that

The deflation-based FastICA functional $\bo\Gamma^d$ satisfies the following $p$ estimating
equations~\cite{Ollila:2010,Nordhausenetal:2011}:
\begin{definition}
\label{defl}
The deflation-based FastICA functional $\bo\Gamma^d=(\bo\gamma_1^d,\dots,\bo\gamma_p^d)^T$ solves the
estimating equations
$$
\bo T(\bo\gamma_j)=\bo S\left(\sum_{l=1}^j\bo\gamma_l\bo\gamma_l^T\right)\,\bo T(\bo\gamma_j),\ \ \
j=1,\dots,p,
$$
where
$$
\bo T(\bo\gamma)=\E[g(\bo\gamma^T(\bo x-\E(\bo x)))(\bo x-\E(\bo x))],
$$
and $g=G'$.
\end{definition}
The estimating equations imply that $\bo \Gamma \bo S\bo \Gamma^T=\bo I_p$, that is, $\bo\Gamma=\bo U\bo S^{-1/2}$
for some orthogonal matrix $\bo U$. The estimation problem can then be reduced to the estimation of the rows of
$\bo U$ one by one. This suggests the following fixed-point algorithm for $\bo u_j$:
\begin{align*}
&\bo u_j\leftarrow \bo T(\bo u_j) \\
&\bo u_j\leftarrow \left(\bo I_p-\sum_{l=1}^{j-1}\bo u_l\bo u_l^T \right)\bo u_j \\
&\bo u_j\leftarrow ||\bo u_j||^{-1}\bo u_j,
\end{align*}
where $\bo T(\bo u)=\E[g(\bo u^T\bo x_{st})\bo x_{st}]$ and  $\bo x_{st}$
is the whitened random variable. However, this algorithm is unstable and we recommend the use of
the original algorithm~\cite{HyvarinenOja:1997}, a modified Newton-Raphson algorithm, where the first step is
$$
\bo u_j\leftarrow \E[g(\bo u_j^T\bo x_{st})\bo x_{st}]-\E[g'(\bo u_j^T\bo x_{st})]\bo u_j.
$$
For the estimate based on the observed data set, all the expectations above are replaced by the sample averages,
e.g., $\E(\bo x)$ is replaced by $\bar{\bo x}$ and $\bo S$ by the sample covariance matrix  $\hat{\bo S}$.

Notice that neither the estimating equations nor the algorithm fixes the order in which the components are
found and the order to some extent depends on the initial value in the algorithm. Since a change in the estimation
order changes the unmixing matrix estimate more than just by permuting its rows, deflation-based FastICA
is not affine equivariant if the initial value is chosen randomly. To find an estimate which
globally maximizes the objective function at each stage, we propose the following strategy to choose the
initial value for the algorithm:
\begin{enumerate}
\item Find a  preliminary consistent estimator $\bo \Gamma_0$ of $\bo \Gamma$.
\item Find a permutation matrix $\bo P$ such that $|\E[G((\bo P\bo \Gamma_0\bo x)_1)]|\geq \cdots \geq |\E[G((\bo P\bo \Gamma_0\bo x)_p)]|$.
\item The orthogonal initial value for $\bo U$ is $\bo P\bo \Gamma_0\bo S^{1/2}$.
\end{enumerate}
The preliminary estimate in step 1 can be for example k-JADE estimate~\cite{Miettinenetal:2013}.
This algorithm, as well as all other FastICA algorithms mentioned in this paper, are implemented in
R package fICA~\cite{Miettinenetal:2014fica}.

The extraction order of the components is highly important not only for the affine equivariance of the estimate,
but also for its efficiency. In the deflationary approach, accurate estimation of the first
components can be shown to have a direct impact on accurate estimation of the last components as well. \cite{Nordhausenetal:2011}
discussed the extraction order and the estimation efficiency and introduced the so-called reloaded deflation-based FastICA,
where the extraction order is based on the minimization of the sum of the asymptotic variances, see
Section~\ref{asymptotics:sec}. \cite{Miettinenetal:2014a} discussed the estimate that uses different G-functions for different
components. Different versions of the algorithm and their performance analysis are presented, for example, in \cite{ZarzasoComonKallel2006,ZarzasoComon2007}.

\subsection{Symmetric approach}
\label{symFastICA:sec}

In symmetric FastICA approach, the rows of $\bo\Gamma=(\bo\gamma_1,\dots,\bo\gamma_p)^T$ are found
simultaneously by maximizing
\[
\sum_{j=1}^p|\E[G(\bo\gamma_j^T(\bo x-\E(\bo x)))]|
\]
under the constraint $\bo\Gamma\bo S\bo\Gamma^T=\bo I_p$.
The unmixing matrix $\bo\Gamma$ optimizes the Lagrangian function
\begin{equation*}
\begin{split}
L(\bo\Gamma,\bo\Theta)=&\sum_{j=1}^p|\E[G(\bo\gamma_j^T(\bo x-\E(\bo x)))]|-\sum_{j=1}^p\theta_{jj}(\bo\gamma^T_j\bo S\bo\gamma_j-1) \\
&-\sum_{j=1}^{p-1}\sum_{l=j+1}^p\theta_{lj}\bo\gamma^T_l\bo S\bo\gamma_j,
\end{split}
\end{equation*}
where symmetric matrix $\bo\Theta=[\theta_{lj}]$ contains $p(p+1)/2$ Lagrangian multipliers. Differentiating
the above function with respect to $\bo\gamma_j$ and setting the derivative to zero yields
\begin{equation*}
\begin{split}
&\ \E[g(\bo\gamma_j^T(\bo x-\E(\bo x)))(\bo x-\E(\bo x))]\,s_j \\ 
=& \ 2\theta_{jj}\bo S\bo\gamma_j
+\sum_{l<j}\theta_{lj}\bo S\bo\gamma_l+\sum_{l>j}\theta_{jl}\bo S\bo\gamma_l,
\end{split}
\end{equation*}
where $g=G'$ and $s_j=sign(\E[G(\bo\gamma_j^T(\bo x-\E(\bo x)))])$. Then by multiplying both
sides by $\bo\gamma_l^T$ we obtain
$\bo\gamma_l^T\E[g(\bo\gamma_j^T(\bo x-\E(\bo x)))(\bo x-\E(\bo x))]s_j=\theta_{lj}$, for $l<j$, and
$\bo\gamma_l^T\E[g(\bo\gamma_j^T(\bo x-\E(\bo x)))(\bo x-\E(\bo x))]s_j=\theta_{jl}$, for $l>j$.
Hence the solution $\bo\Gamma$ must satisfy the following estimating equations
\begin{definition}
\label{def2}
The symmetric FastICA functional $\bo\Gamma^s=(\bo\gamma_1^s,\dots,\bo\gamma_p^s)^T$ solves the estimating
equations
$$
\bo\gamma_l^T\bo T(\bo\gamma_j)\,s_j=\bo\gamma_j^T\bo T(\bo\gamma_l)\,s_l\ \ \
\text{and}\ \ \ \bo\gamma^T_lS\bo\gamma_j=\delta_{lj},
$$
where $j,l=1,\dots,p$, and
$$
\bo T(\bo\gamma)=\E[g(\bo\gamma^T(\bo x-\E(\bo x)))(\bo x-\E(\bo x))],
$$
$g=G'$, $s_j=sign(\E[G(\bo\gamma_j^T(\bo x-\E(\bo x)))]$ and $\delta_{lj}$ is the Kronecker delta.
\end{definition}
Again, $\bo\Gamma=\bo U\bo S^{-1/2}$ for some orthogonal matrix $\bo U$.
Then the estimation equations for $\bo U$ are
\[
 \bo u_l^T\bo T(\bo u_j)\,s_j=\bo u_j^T\bo T(\bo u_l)\,s_l\ \ \ \text{and}\ \ \
\bo u^T_l\bo u_j=\delta_{lj},
\]
where $l,j=1,\dots,p$, $\bo T(\bo u)=\E[g(\bo u^T\bo x_{st})\bo x_{st}]$, and the equations suggest the following fixed-point algorithm
for $\bo U$:
 \begin{align*}
& \bo T\leftarrow (\bo T(\bo u_1),\dots,\bo T(\bo u_p))^T \\
& \bo U\leftarrow (\bo T\bo T^T)^{-1/2} \bo T.
\end{align*}
As in the deflation-based approach, a more stable algorithm is obtained when $\bo T(\bo u_j)$
is replaced by
$$
\bo T^*(\bo u_j)=\E[g(\bo u_j^T\bo x_{st})\bo x_{st}]-\E[g'(\bo u_j^T\bo x_{st})]\bo u_j.
$$
In symmetric FastICA, different initial values give identical unmixing matrix estimates up to order
and signs of the rows.

\subsection{Squared symmetric approach}
\label{squaredsymFastICA:sec}

In squared symmetric FastICA, the absolute values in the objective function of the regular
symmetric FastICA are replaced by squares~\cite{VirtaNordhausenOja:2015}. The squared symmetric FastICA functional
$\bo\Gamma^{s2}=(\bo\gamma_1^{s2},\dots,\bo\gamma_p^{s2})^T$ maximizes
\[
\sum_{j=1}^p(\E[G(\bo\gamma_j^T(\bo x-\E(\bo x)))])^2
\]
under the constraint $\bo\Gamma\bo S\bo\Gamma^T=\bo I_p$.
Similarly as in Section~\ref{symFastICA:sec}
the Lagrange multipliers method yields the following estimating equations:
\begin{definition}
\label{def3}
The squared symmetric FastICA functional $\bo\Gamma^{s2}=(\bo\gamma_1^{s2},\dots,\bo\gamma_p^{s2})^T$ solves
the estimating equations
$$
\bo\gamma_l^T\bo T_2(\bo\gamma_j)=\bo\gamma_j^T\bo T_2(\bo\gamma_l)\ \ \
\text{and}\ \ \ \bo\gamma^T_l\bo S\bo\gamma_j=\delta_{lj},
$$
where $j,l=1,\dots,p$,
$$
\bo T_2(\bo\gamma)=\E[G(\bo\gamma^T(\bo x-\E(\bo x)))]\E[g(\bo\gamma^T(\bo x-\E(\bo x)))(\bo x-\E(\bo x))],
$$
 $g=G'$ and $\delta_{lj}$ is the Kronecker delta.
\end{definition}
The estimation equations for $\bo U$ are
\[
 \bo u_l^T\bo T_2(\bo u_j)=\bo u_j^T\bo T_2(\bo u_l)\ \ \ \text{and}\ \ \
\bo u^T_l\bo u_j=\delta_{lj}, \ \ l,j=1,\dots,p,
\]
where $\bo T_2(\bo u)=\E[G(\bo u^T(\bo x_{st}))]\E[g(\bo u^T(\bo x_{st}))\bo x_{st}]$.
The following algorithm, which is based on the same idea as the algorithm for symmetric FastICA,
can be used to find the solution in practice:
\begin{align*}
& \bo T\leftarrow (\bo T_2^*(\bo u_1),\dots,\bo T_2^*(\bo u_p))^T \\
& \bo U\leftarrow (\bo T\bo T^T)^{-1/2} \bo T,
\end{align*}
where $\bo T_2^*(\bo u)=\E[G(\bo u^T(\bo x_{st}))]\{\E[g(\bo u^T(\bo x_{st}))\bo x_{st}]-\E[g'(\bo u^T(\bo x_{st}))]\bo u \}$.

Notice that
\[
\bo T_2^*(\bo u)=\E[G(\bo u^T(\bo x_{st}))]\bo T^*(\bo u),
\]
and hence the squared symmetric FastICA estimator can be seen as weighted classical symmetric FastICA estimator. The more
nongaussian, as measured by function $G$, an independent component is, the more impact it has in the orthogonalization step.

\subsection{Function $G$}
\label{Gfunction:sec}

The function $G$  is required to be twice continuously differentiable, nonlinear and nonquadratic
function such that $\E[G(z)]=0$, when $z$ is a standard Gaussian random variable.
The derivative function $g=G'$ is the  so-called nonlinearity. The use of classical kurtosis as a measure of non-Gaussianity
is given by the nonlinearity function $g(z)=z^3$ ({\it pow3})~\cite{HyvarinenOja:1997}. Other popular choices
include $g(z)=tanh(az)$ ({\it tanh}) and $g(z)=z\,\exp(-az^2/2)$ ({\it gaus}) with tuning parameters
$a$ as suggested in~\cite{Hyvarinen:1997}, and $g(z)=z^2$ ({\it skew}).

\begin{figure}[htb]
\centering
\includegraphics[width=1.7in]{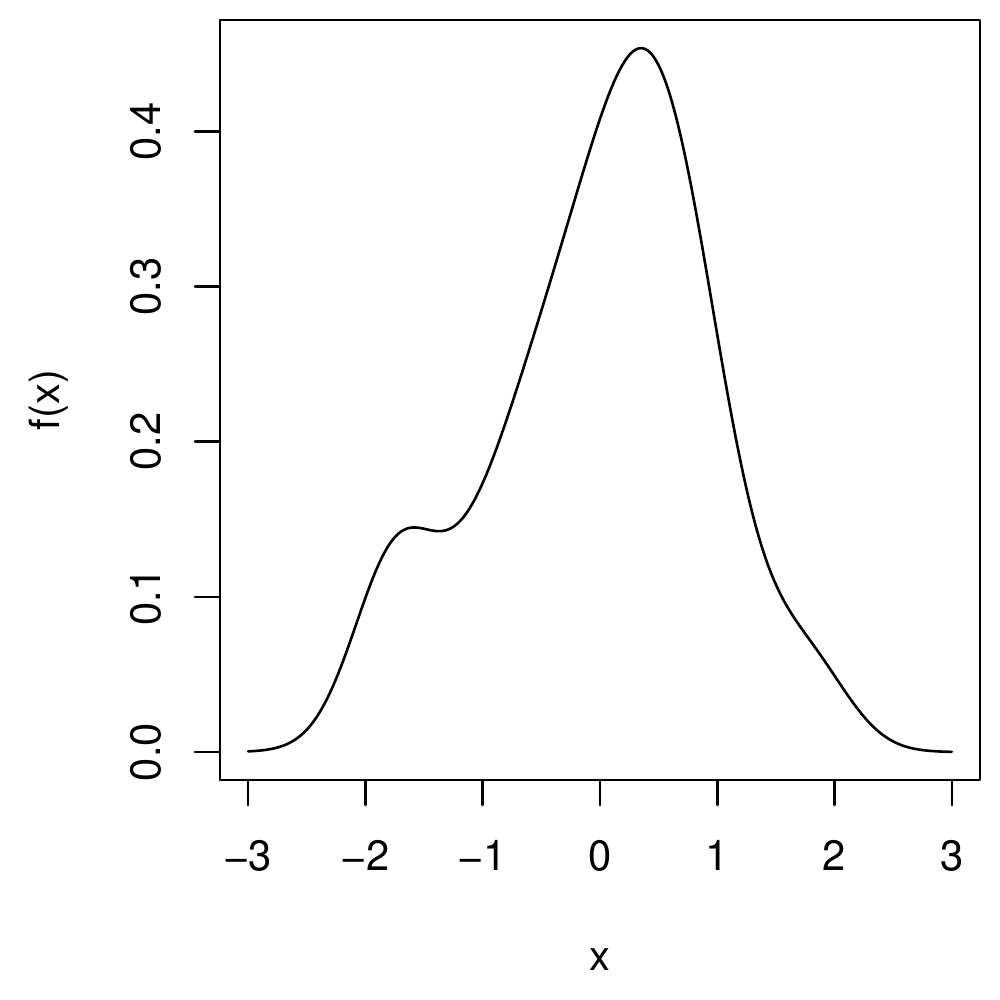}
\includegraphics[width=1.7in]{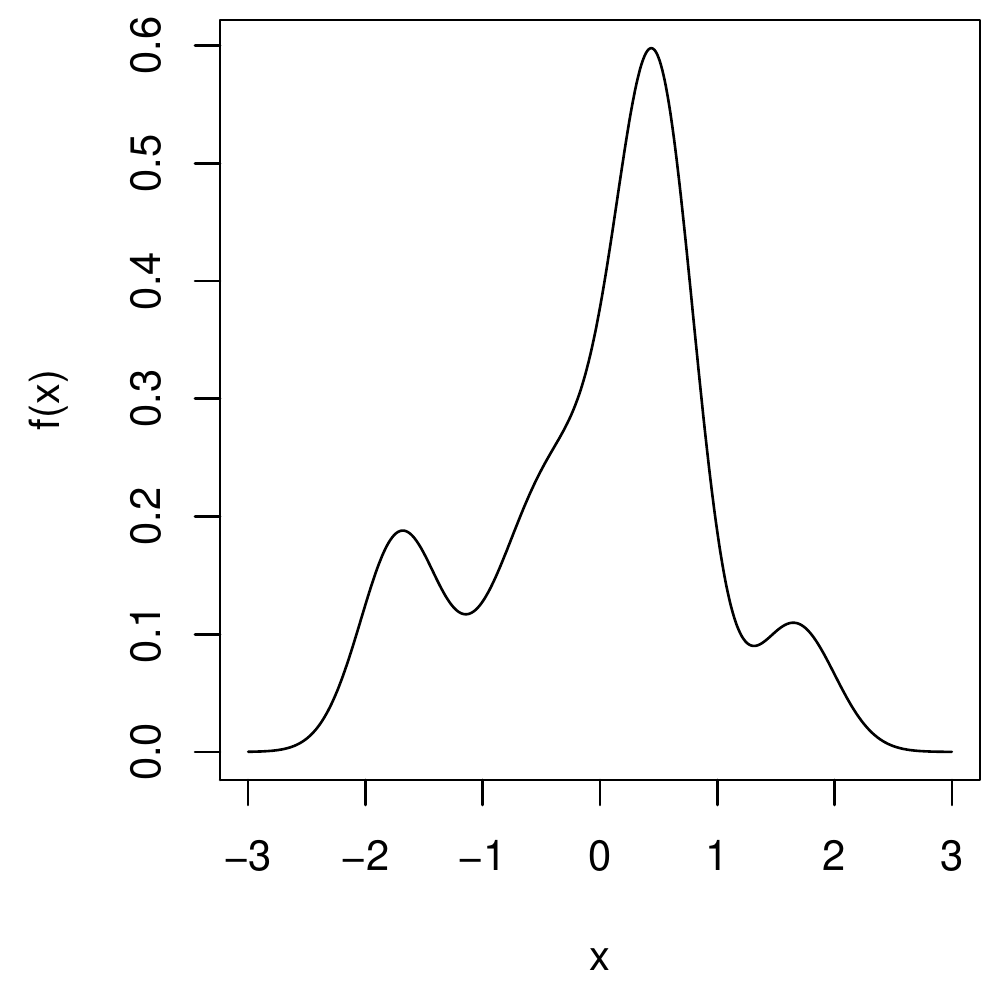}
\caption{Density functions of $z_1$ and $z_2$, which violate the conditions def, sym and sym2 with nonlinearity $gaus$. Both distributions are mixtures of four Gaussian distributions. For more details, see Appendix. }
\label{fig1}
\end{figure}

The deflation-based, symmetric and squared symmetric FastICA estimators need extra conditions for $G$  to ensure the consistency of the estimation procedure:
One then requires that, for any bivariate  $\bo Z=(z_1,z_2)^T$  with independent and standardized components ($\E(\bo z)=\bo 0$ and $\cov(\bo z)=\bo I_2$)  and for any orthogonal $2\times 2$ matrix $\bo U=(\bo u_1,\bo u_2)^T$,
\begin{align*}
&\mbox{def} \ \ \ \ |\E[G(\bo u_{1}^T \bo z)]|\leq \mbox{max}(|\E[G(z_1)]|,|\E[G(z_2)]|), \\
&\mbox{sym} \ \ \ |\E[G(\bo u_1^T\bo z)]|+ |\E[G(\bo u_2^T\bo z)]| \\
&\quad\quad\quad  \leq |\E[G(z_1)]|+ |\E[G(z_2)]| \\
&\mbox{sym2} \ \ \ \left(\E[G(\bo u_1^T\bo z)]\right)^2+ \left(\E[G(\bo u_2^T \bo z)]\right)^2\\
&\quad\quad\quad    \leq \left(\E[G(z_1)]\right)^2+ \left(\E[G(z_2)]\right)^2
\end{align*}
\cite{Miettinenetal:2015} and \cite{VirtaNordhausenOja:2015}
proved that for {\it pow3} and  {\it skew} (as well as for their convex combination),  all three conditions  are satisfied.
On the contrary, {\it tanh} and {\it gaus} do not satisfy the conditions for all choices of the distributions of $z_1$ and $z_2$. For these two
nonlinearities \cite{Wei:2014a} found bimodal distributions for which the fixed points of the deflation-based
FastICA algorithm are not correct solutions of the IC problem. In Figure~\ref{fig1} we plot the density
functions of random variables $z_1$ and $z_2$ which serve as examples for a case where none of the three
inequalities hold for {\it gaus}.
These examples should however be seen as rare and artificial exceptions and FastICA with {\it tanh} and {\it gaus}
 satisfy the conditions for most of pairs of distributions of $\bo z$ we have checked. For example, in Section~\ref{simulation:sec}
FastICA with {\it tanh} worked as expected under a wide variety of source distributions.
Deflation-based or symmetric FastICA with {\it tanh} is perhaps the most popular unmixing matrix estimate.

See Section~\ref{asympeff:sec} for the optimal choice of the nonlinearity for a component with a known
density function.

\section{Asymptotical properties of the FastICA estimators}
\label{asymptotics:sec}

The limiting variances and the asymptotic multinormality of the deflation-based and symmetric FastICA unmixing matrix
estimators were found quite recently in~\cite{Ollila:2010},~\cite{Nordhausenetal:2011},~\cite{Wei:2014b} and~\cite{Wei:2015}.
In this section, we review these findings and derive the results for the squared symmetric
FastICA estimator.

Let now $\bo X=(\bo x_1,\dots,\bo x_n)$ be a random sample from the distribution of $\bo x$ following
the IC model~(\ref{ICmodel}). The deflation-based, symmetric and squared symmetric FastICA estimators
$\hat{\bo\Gamma}^d$, $\hat{\bo\Gamma}^s$ and $\hat{\bo\Gamma}^{s2}$
are then obtained when the three functionals are applied to the empirical distribution of $\bo X$.

Due to affine equivariance, we can in the following assume without loss of generality that
$\bo\Omega=\bo I_p$. Before proceeding we need to make some additional assumptions on the distribution of
$\bo z_i=(z_{i1},\dots,z_{ip})^T$, namely,
\begin{itemize}
\item[(A3)] The fourth moments $\beta_j=\E[z_{ij}^4]$ as well as the following expected values
\begin{align*}
\nu_j&=\E[G(z_{ij})], &\mu_j&=\E[g(z_{ij})], &\sigma_j^2&=\mbox{Var}[g(z_{ij})], \\
\lambda_j&=\E[g(z_{ij})z_{ij}], &\delta_j&=\E[g'(z_{ij})], &\tau_j&=\E[g'(z_{ij})z_{ij}]
\end{align*}
exist. Write also $s_j=sign(\nu_j)$.
\end{itemize}

Write now
\begin{align*}
\bo T_j&=\frac{1}{n}\sum_{i=1}^n(g(z_{ij})-\mu_j)\bo z_i \ \ \ \ \text{ and} \\
\bo T_{2j}&=\frac{1}{n}\sum_{i=1}^nG(z_{ij})\frac{1}{n}\sum_{i=1}^n(g(z_{ij})-\mu_j)\bo z_i
\end{align*}
for $j=1,\dots,p$. To avoid division by zero in the following theorem, assume that $\nu_j(\lambda_j-\delta_j)\geq 0$ for all
$j=1,\dots,p$, with equality for at most one $j$, see \cite{Hyvarinen:1999}, who stated that
$\nu_j(\lambda_j-\delta_j)>0$ for most of the reasonable functions $G$ and distributions of $z_{ij}$.
For $(pow3)$, $\nu_j(\lambda_j-\delta_j)>0$ for any distribution with $\E(z_{ij}^4)\neq 3$.
The limiting behavior of the deflation-based FastICA estimate was first given in~\cite{Nordhausenetal:2011}.
The corresponding results of the symmetrical FastICA estimates are given in the following. The result (iii)
is proved in the Appendix and the proof of (ii) is essentially similar to that.  In the
following theorem,  $\bo e_i$ is a $p$-vector with $i$th element one and others zero and $o_P(1)$ replaces a random variable that converges in probability to zero
as $n$ goes to the infinity.
\begin{theorem}
\label{asymp_theorem}
Let $\bo X=(\bo x_1,\dots,\bo x_n)$ be a random sample from the IC model~(\ref{ICmodel}) satisfying
the assumptions (A1)-(A3): If $\bo \Omega=\bo I_p$ then there exist a sequence of solutions $\hat{\bo\Gamma}^d$,
$\hat{\bo\Gamma}^s$ and $\hat{\bo\Gamma}^{s2}$ converging to $\bo I_p$ such that

(i) (deflation-based)
\begin{align*}
\begin{split}
&\sqrt{n}\,{\hat\gamma}_{jl}^d = -\sqrt{n}\,{\hat\gamma}_{lj}^d-\sqrt{n}\, \hat{\bo S}_{jl}+o_P(1),\ \ \ \ l<j, \\
&\sqrt{n}\,({\hat\gamma}_{jj}^d-1) = -\frac{1}{2}\sqrt{n}\,(\hat{\bo S}_{jj}-1)+o_P(1),\ \ \ \ l=j, \\
&\sqrt{n}\,{\hat\gamma}_{jl}^d = \frac {\bo e_l^T \sqrt{n}\,\bo T_{j}-\lambda_j\sqrt{n}\,\hat{\bo S}_{jl}}
{\lambda_j-\delta_j} +o_P(1),\ \ \ \ l>j,
\end{split}
\end{align*}

(ii) (symmetric)
\begin{equation*}
\begin{split}
&\sqrt{n}\,(\hat{\gamma}_{jj}^s-1)=-\frac{1}{2}\sqrt{n}\,(\hat{\bo S}_{jj}-1)+o_P(1),\ \ \ \ l=j, \\  \\
&\sqrt{n}\,\hat{\gamma}_{jl}^s = \frac{\bo e_l^T\sqrt{n}\,\bo T_js_j-\bo e_j^T\sqrt{n}\,\bo T_ls_l-(\lambda_j\,s_j-\delta_ls_l)\sqrt{n}\,\hat{\bo S}_{jl}}{(\lambda_j-\delta_j)s_j+(\lambda_l-\delta_l)s_l} \\ &\qquad\qquad+o_P(1),\ \ \ \ l\neq j,
\end{split}
\end{equation*}

(iii) (squared symmetric)
\begin{equation*}
\begin{split}
&\sqrt{n}\,(\hat{\gamma}_{jj}^{s2}-1)=-\frac{1}{2}\sqrt{n}\,(\hat{\bo S}_{jj}-1)+o_P(1),\ \ \ \ l=j, \\  \\
&\sqrt{n}\,\hat{\gamma}_{jl}^{s2} = \frac{\bo e_l^T\sqrt{n}\,\bo T_{2j}-\bo e_j^T\sqrt{n}\,\bo T_{2l}+
(\nu_l\delta_l-\nu_j\lambda_j)\sqrt{n}\,\hat{\bo S}_{jl}}{\nu_j(\lambda_j-
\delta_j)+\nu_l(\lambda_l-\delta_l)} \\ &\qquad\qquad+o_P(1),\ \ \ \ l\neq j.
\end{split}
\end{equation*}
\end{theorem}
For the asymptotical properties of deflation-based FastICA for several nonlinearities $g$,
see~\cite{Miettinenetal:2014a}.
As seen from Theorem~\ref{asymp_theorem} (i), the limiting distributions of vectors $\hat{\bo\gamma}_1^d,\dots,\hat{\bo\gamma}_p^d$
depend on the order in which they are found. It is shown in 
Corollary~\ref{asv_cor} that, for $j<l$,  the asymptotic variances of $\hat{\gamma}_{lj}^d$ and $\hat{\gamma}_{jl}^d$ are equal and 
depend only on the distribution of the $j$th independent component. The limiting distributions of the diagonal elements
do not depend on the method or the chosen nonlinearity $g$.
%, but they are induced by the constraint $\bo \Gamma\bo S\bo \Gamma^T=\bo I_p$.
\cite{VirtaNordhausenOja:2015} discovered that the squared symmetric FastICA estimator with $(pow3)$ nonlinearity has the same asymptotics as the  JADE (joint approximate diagonalization of eigenmatrices) estimator~\cite{CardosoSouloumiac:1993}. We then have the following straightforward but important corollaries.
\begin{corollary}
\label{asnormal}
Under the assumptions of Theorem~\ref{asymp_theorem}, if the joint limiting distribution of $\sqrt{n}\,\bo T_{jl}$
and  $\sqrt{n}\,\bo T_{2jl}$ for $j\ne l=1,\dots,p$ and $\sqrt{n}\,(\hat{\bo S}_{jl}-\delta_{jl})$ for $j,l=1,\dots,p$, is a
 multivariate normal distribution, then also the limiting distributions of
$\sqrt{n}\,vec(\hat{\bo\Gamma}^d-\bo I_p)$,  $\sqrt{n}\,vec(\hat{\bo\Gamma}^s-\bo I_p)$ and $\sqrt{n}\,vec(\hat{\bo\Gamma}^{s2}-\bo I_p)$ are multivariate normal.
\end{corollary}

\begin{corollary}
\label{asv_cor}
Under the assumptions of Theorem~\ref{asymp_theorem}, the asymptotic covariance matrix (ASV)
of the $j$th source vectors are given by
\begin{align*}
ASV(\hat{\bo\gamma}_{j}^d)&=\sum_{l=1}^p ASV(\hat{\gamma}_{jl}^d)\bo e_l\bo e_l^T, \\
ASV(\hat{\bo\gamma}_{j}^s)&=\sum_{l=1}^p ASV(\hat{\gamma}_{jl}^s)\bo e_l\bo e_l^T, \ \ \text{ and} \\
ASV(\hat{\bo\gamma}_{j}^{s2})&=\sum_{l=1}^p ASV(\hat{\gamma}_{jl}^{s2})\bo e_l\bo e_l^T,
\end{align*}
where

(i) (deflation-based)
\begin{align*}
\begin{split}
& ASV(\hat{\gamma}_{jl}^d)=\frac{\sigma^2_l-\lambda^2_l}{(\lambda_l-\delta_l)^2}+1,\ \ \ \ l<j \\
& ASV(\hat{\gamma}_{jj}^d)=\frac{\beta_j-1}{4},\ \ \ \ l=j \\
& ASV(\hat{\gamma}_{jl}^d)=\frac{\sigma^2_j-\lambda^2_j}{(\lambda_j-\delta_j)^2},\ \ \ \ l>j.
\end{split}
\end{align*}

(ii) (symmetric)
\begin{align*}
\begin{split}
& ASV(\hat{\gamma}_{jj}^s)=\frac{\beta_j-1}{4},\ \ \ \ l=j \\
& ASV(\hat{\gamma}_{jl}^s)=\frac{\sigma_j^2+\sigma_l^2-\lambda^2_j+\delta_l(\delta_l-2\lambda_l)}{((\lambda_j-\delta_j)s_j+(\lambda_l-\delta_l)s_l)^2},\ \ \ \ l\neq j.
\end{split}
\end{align*}

(iii) (squared symmetric)
\begin{align*}
\begin{split}
 ASV(\hat{\gamma}_{jj}^{s2})=&\frac{\beta_j-1}{4},\ \ \ \ l=j \\
 ASV(\hat{\gamma}_{jl}^{s2})=&\frac{\nu_j^2(\sigma_j^2-\lambda^2_j)+\nu_l^2(\sigma_l^2+
\delta_l(\delta_l-2\lambda_l))}{(\nu_j(\lambda_j-\delta_j)+
\nu_l(\lambda_l-\delta_l))^2},\\
& l\neq j.
\end{split}
\end{align*}
\end{corollary}
The asymptotic variances of the deflation-based and symmetric FastICA estimators were first derived
in \cite{Ollila:2010} and \cite{Wei:2014b}, respectively. The asymptotic covariance matrices of the
FastICA estimators for given marginal densities can be computed  using the R package BSSasymp \cite{Miettinenetal:2014BSSsymp}.

\section{Efficiency comparisons}
\label{comparisons:sec}

The asymptotical results derived in Section~\ref{asymptotics:sec} allow us to evaluate and compare the
performances of the FastICA methods. In this section the asymptotic and finite sample efficiencies of
deflation-based and symmetric FastICA  estimators are compared to those of squared symmetric FastICA
estimators using a wide range of distributions with varying skewness and kurtosis values.

\subsection{Performance index}
\label{mdsec}
We measure the finite sample performance of the unmixing matrix estimates using the minimum distance
index~\cite{IlmonenNordhausenOjaOllila:2010b}
\begin{align}
\label{md}
\hat{D}= D(\hat{\bo\Gamma}\bo\Omega)=\frac{1}{\sqrt{p-1}}\inf_{\bo C\in
\mathcal{C}}\|\bo C\hat{\bo\Gamma}\bo\Omega-\bo I_p\|
\end{align}
where $\|\cdot\|$ is the matrix (Frobenius) norm and $\mathcal{C}$ is the set of $p\times p$ matrices with exactly one
non-zero element in each column and each row. The minimum distance index is scaled so that $0\leq\hat{D}\leq 1$.
If $\bo\Omega=\bo I_p$ and $\sqrt{n}\,vec(\hat{\bo\Gamma}-\bo I_p)\to N_{p^2}(\bo 0,\bo\Sigma)$, then the limiting
distribution of $n(p-1)\hat{D}^2$ is that of weighted sum of independent chi squared variables with the expected value
\begin{align}
\label{md2}
Trace\left[(\bo I_{p^2}-\bo D_{p,p})\bo\Sigma (\bo I_{p^2}-\bo D_{p,p})\right],
\end{align}
where $\bo D_{p,p}=\sum_i(\bo e_i\bo e_i^T)\otimes(\bo e_i\bo e_i^T)$, and $\otimes$ means the Kronecker product.
Notice that~(\ref{md2}) equals the sum of the limiting variances
of the off-diagonal elements of $\sqrt{n}\,{vec}(\hat{\bo\Gamma}-\bo I_p)$ and therefore
\begin{align}
\label{md3}
\sum_{j=1}^{p-1}\sum_{l=j+1}^p \left(ASV(\hat\gamma_{jl})+ASV(\hat\gamma_{lj})\right)
\end{align}
provides a global measure of the variation of the estimate $\hat{\bo\Gamma}$.

\subsection{Asymptotic efficiency}
\label{asympeff:sec}

Let $f_j$ be the density function and  $g_j=-f'_j/f_j$ be the optimal location score function for the $j$th independent
component $z_j$. Also let $I_j=Var(g_j(z_j))$ be the Fisher information number for the location problem. Write
\[
\alpha_j:=\frac{\sigma_j^2-\lambda_j^2}{(\lambda_j-\delta_j)^2}=[(I_j-1) \rho^2_{g(z_j)g_j(z_j)\cdot z_j} ]^{-1},
\]
where $\rho^2_{g(z_j)g_j(z)\cdot z_j}$ is the squared partial correlation between $g(z_j)$ and $g_j(z_j)$ given $z_j$.
Then we have the following.
\begin{theorem}
For our three estimates and for non-gaussian $z_j$ and $z_l$, $j\ne l$,
$ASV(\hat\gamma_{jl})+ASV(\hat\gamma_{lj})$ is
\[
\left(\frac {\beta_j}{\beta_j+\beta_l}\right)^2 (2\alpha_j+1)+
\left(\frac {\beta_l}{\beta_j+\beta_l}\right)^2 (2\alpha_l+1)
\]
where
\[
\left\{ \begin{array}{ll}
 \beta_j &=1,\ \  s_j(\lambda_j-\delta_j),\ \mbox{and}\ \  \nu_j(\lambda_j-\delta_j)   \\
 \beta_l &=0,\ \ s_l(\lambda_l-\delta_l),\ \mbox{and}\ \ \nu_l(\lambda_l-\delta_l)  \end{array} \right.
\]
for deflation-based, symmetric and squared symmetric FastICA estimates, respectively.
\end{theorem}

Notice first that the value of $ASV(\hat\gamma_{jl})+ASV(\hat\gamma_{lj})$ only depends on the $j$th and $l$th marginal distributions, which means we can restrict the comparison to bivariate distributions as the other components have no impact. If the $j$th and $l$th marginal distributions are the same, then the three values of  $ASV(\hat\gamma_{jl})+ASV(\hat\gamma_{lj})$ are
\[
(2\alpha_j+1), \ \ \frac 12 (2\alpha_j+1)\ \ \mbox{and}\ \ \frac 12 (2\alpha_j+1)
\]
and these are minimized with the choice $g=g_j$. So, if $z_1,\dots,z_p$ are identically distributed with the density function $f$, then the optimal choice for $g$ is $-f'/f$.

If the $l$th component is  Gaussian then, $\lambda_l=\delta_l$, and for the deflation-based and squared symmetric FastICA
estimates,  $ASV(\hat\gamma_{jl})+ASV(\hat\gamma_{lj})=(2\alpha_j+1)$ and for the symmetric FastICA estimate one gets
\begin{eqnarray*}
&&ASV(\hat\gamma_{jl})+ASV(\hat\gamma_{lj})
=(2\alpha_j+1)+\frac{2(\sigma_l^2-\lambda_l^2)}{\beta_j^2}\\
&=&(2\alpha_j+1)+\frac{2\sigma_l^2}{\beta_j^2} \left(1 -\rho_{g(z_{il})z_{il}}^2\right)
\end{eqnarray*}
where $\rho_{g(z_{il})z_{il}}$ is the correlation between $g(z_{il})$ and $z_{il}$. The symmetric FastICA is therefore always poorest in this case.

For further comparison of the estimators we use two families of source distributions, the
standardized exponential power distribution family and the standardized gamma distribution family.
The density function of standardized exponential power distribution with shape parameter $\beta$ is
\[
f(x)=\frac{\beta\,\exp\{-(|x|/\alpha)^\beta\}}{2\,\alpha\,\Gamma(1/\beta)},
\]
where $\beta>0$, $\alpha=(\Gamma(1/\beta)/\Gamma(3/\beta))^{1/2}$ and $\Gamma$ is the gamma function.
The distribution is symmetric for any $\beta$, and $\beta=2$ gives the normal (Gaussian) distribution,
$\beta=1$ gives the heavy-tailed Laplace distribution and the density converges to the low-tailed uniform
distribution as $\beta\to\infty$.
The density function of standardized gamma distribution with shape parameter $\alpha$ is
\[
f(x)=\frac{(x+\sqrt{\alpha})^{\alpha-1}\alpha^{\alpha/2} \exp\{-(x+\sqrt{\alpha})\sqrt{\alpha}\}}{\Gamma(\alpha)}.
\]
Gamma distributions are right skew, and for $\alpha=k/2$, the distribution is a chi square distribution
with $k$ degrees of freedom, $k=1,2,\dots$. When $\alpha=1$, we have an exponential distribution, and
the distribution converges to a normal distribution as $\alpha\to\infty$.

We next compare the asymptotic variances of the unmixing matrix estimates with the same nonlinearity and for $\bo\Omega=\bo I_p$. For the comparison, write
\[
ARE_{s2,d}=\frac{ASV(\hat{\gamma}_{jl}^d)+ASV(\hat{\gamma}_{lj}^d)}{ASV(\hat{\gamma}_{jl}^{s2})+ASV(\hat{\gamma}_{lj}^{s2})},
\]
for the asymptotic relative efficiency of the squared symmetric estimate with respect to the deflation based estimate, and similarly
for $ARE_{s2,s}$. Notice that $ARE_{s2,d}$ and $ARE_{s2,s}$ depend on the two marginal distribution as well as on the chosen nonlinearity. 
We then  plot the contour maps of the ARE's as functions of the shape parameters of the  exponential power or gamma distributions with nonlinearities
 {\it pow3} and {\it tanh}.  The equal efficiency is given by the ARE value 1 and can be found using the bar with contour thresholds on the
right-hand side of the figures.

\begin{figure*}[!t]
\centering
\includegraphics[width=8.5in]{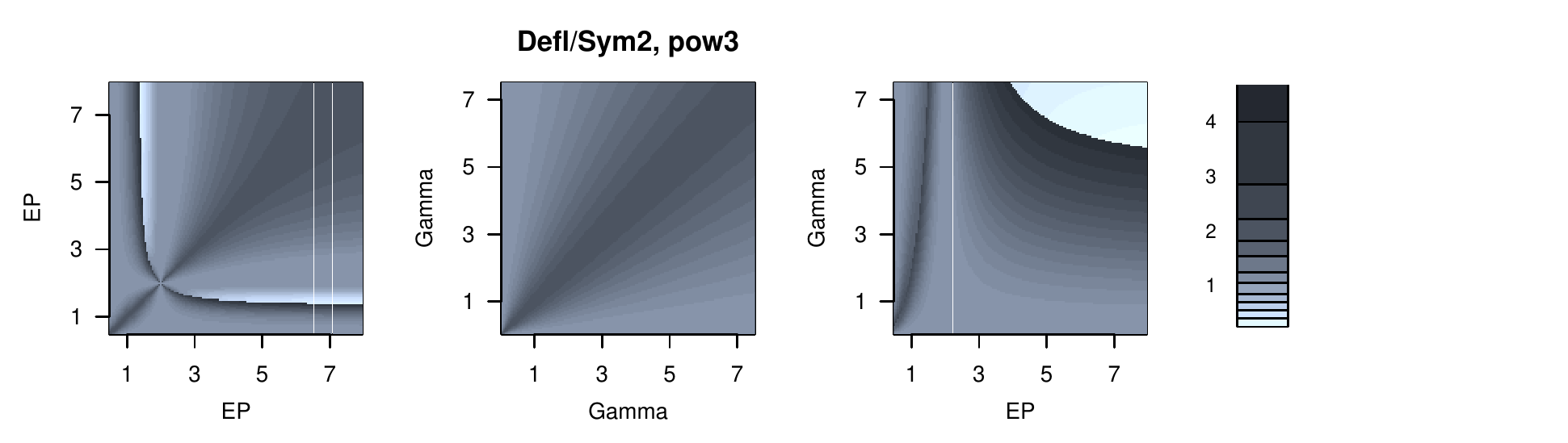}

\includegraphics[width=8.5in]{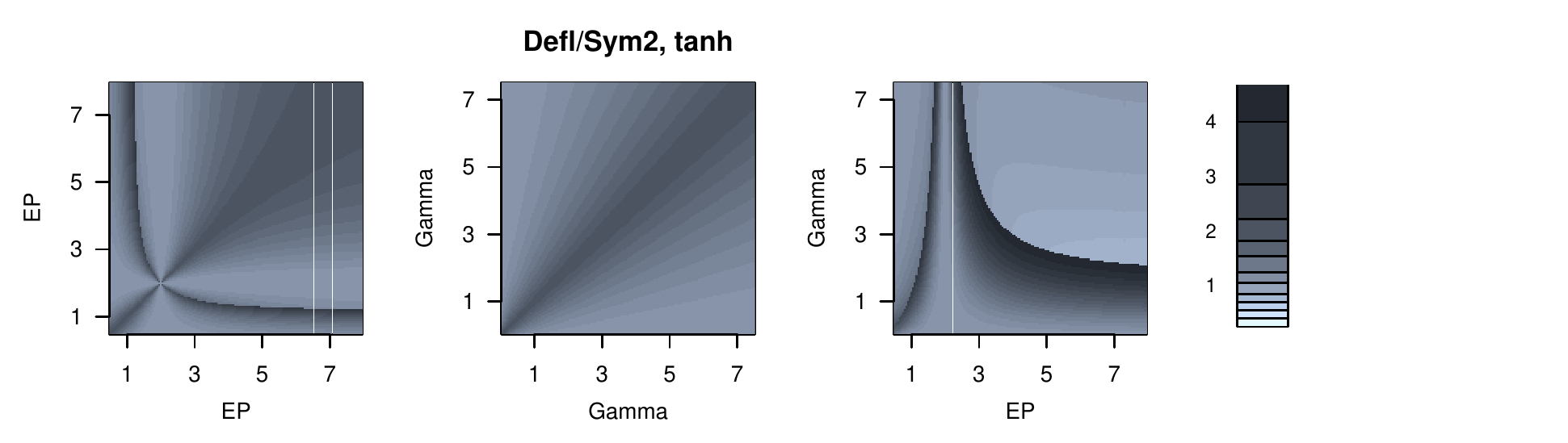}
\caption{Contour maps of asymptotic relative efficiencies  $ARE_{s2,d}$   when $\bo\Omega=\bo I_p$ and
the source distributions are exponential power (EP) or gamma (Gamma) distributions with varying
shape parameter values. The nonlinearities are \textit{pow3} on the top row and \textit{tanh} on the bottom row.}
\label{fig2}
\end{figure*}

\begin{figure*}[!t]
\centering
\includegraphics[width=8.5in]{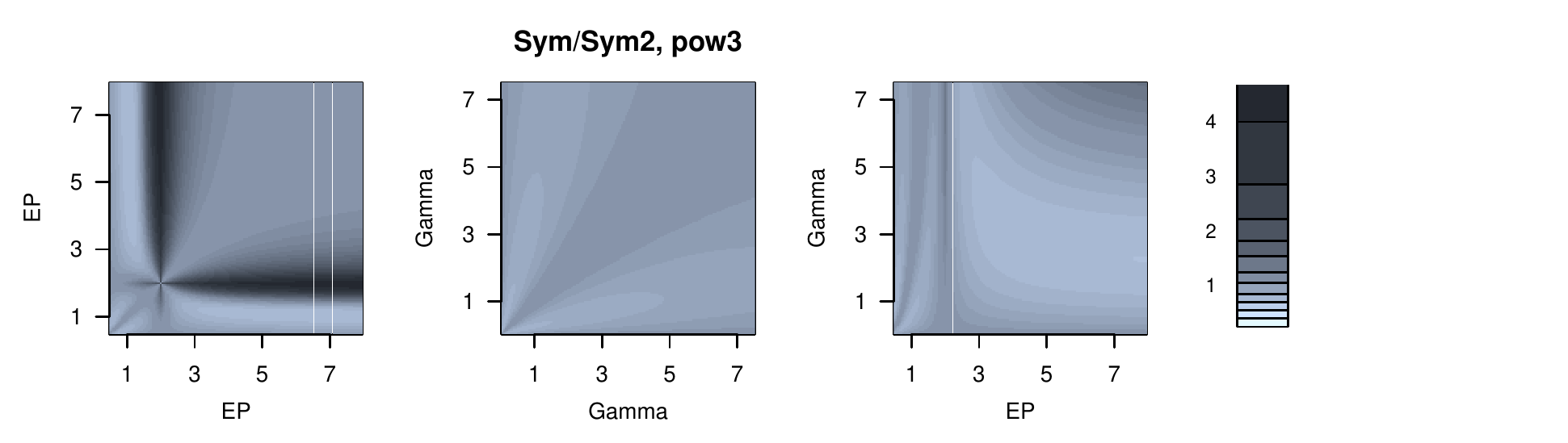}

\includegraphics[width=8.5in]{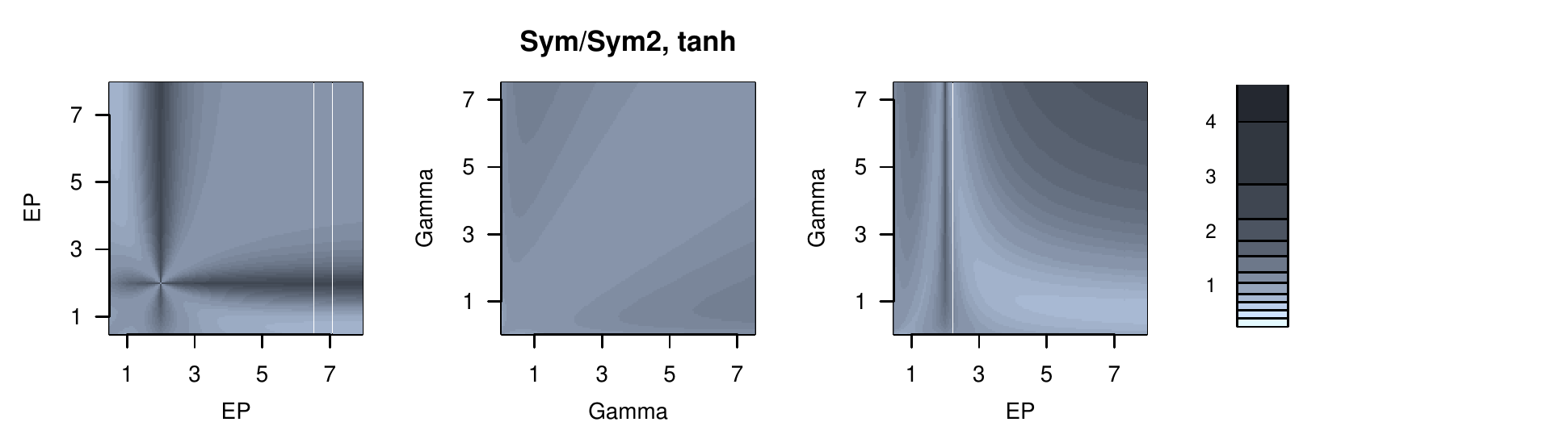}
\caption{Contour maps of asymptotic relative efficiencies  $ARE_{s2,s}$   when $\bo\Omega=\bo I_p$ and
the source distributions are exponential power (EP) or gamma (Gamma) distributions with varying
shape parameter values. The nonlinearities are \textit{pow3} on the top row and \textit{tanh} on the bottom row.}
\label{fig3}
\end{figure*}

As seen in Figure~\ref{fig2}, the squared symmetric FastICA estimator is in most cases more efficient than the deflation-based estimator. In Figure~\ref{fig3} we use $ARE_{s2,s}$ similarly for the comparison between symmetric and squared symmetric FastICA. In the figures, the darker the point the higher relative efficiency. Notice that $ARE_{s2,d}=1$  if one of the components is Gaussian, and  $ARE_{s2,s}=1$ if $E(G(z_1))=E(G(z_2))$ (e.g. if the two distributions are the same).
Figure~\ref{fig3} shows that the areas where
$ASV(\hat{\gamma}_{jl}^s)>ASV(\hat{\gamma}_{jl}^{s2})$ and $ASV(\hat{\gamma}_{jl}^s)<ASV(\hat{\gamma}_{jl}^{s2})$ are almost equally large, but the differences in favour of the squared symmetric estimator are larger. Also, they occur in cases where the separation of the components is difficult, and hence the efficiency is important there.

In Table~\ref{table1} the values of  $ARE_{s2,s}$  and  $ARE_{s2,d}$  are displayed for different pairs of source distributions and for  \emph{pow3} in the upper triangle and for \emph{tanh} in the lower triangle. Table~\ref{table1} presents a sample of the values of Figure~\ref{fig2} and Figure~\ref{fig3} in a numerical form.

\begin{table}[htb]
\renewcommand{\arraystretch}{1.3}
\caption{Values of $ARE_{s2,s}$  (on the top) and  $ARE_{s2,d}$ (on the bottom) for different distributions. L=Laplace=EP1=exponential power distribution with $\beta=1$, N=Normal distribution=EP2, U=Uniform distribution, G1=Gamma distribution with $\alpha=1$. Upper triangle for $pow3$ and lower triangle for $tanh$.}
\label{table1}
\centering
  \begin{tabular*}{\columnwidth}{@{}*{11}{c@{\extracolsep\fill}}c@{}}
    \hline
        & L    & EP1.5 & EP1.75 & N    & EP 3 & EP4  & U    & G1   & G3   & G6   \\
 L      & 1    & 0.87  & 0.96   & 1.10 & 0.79 & 0.74 & 0.70 & 0.84 & 1.00 & 0.94 \\
 EP1.5  & 1.05 & 1     & 1.02   & 1.46 & 0.87 & 1.05 & 1.61 & 0.87 & 0.84 & 0.95 \\
 EP1.75 & 1.21 & 1.10  & 1      & 1.72 & 1.53 & 2.21 & 3.75 & 0.95 & 0.93 & 0.92 \\
 N      & 1.50 & 1.79  & 1.90   & --   & 3.09 & 3.85 & 5.49 & 1.03 & 1.13 & 1.25 \\
 EP3    & 0.91 & 0.99  & 1.29   & 2.23 & 1    & 1.23 & 1.45 & 0.86 & 0.74 & 0.77 \\
 EP4    & 0.85 & 1.13  & 1.55   & 2.32 & 1.05 & 1    & 1.13 & 0.81 & 0.71 & 0.86 \\
 U      & 0.87 & 1.46  & 1.90   & 2.45 & 1.24 & 1.08 & 1    & 0.76 & 0.74 & 1.20 \\
 G1     & 0.97 & 0.94  & 1.12   & 1.43 & 0.82 & 0.74 & 0.76 & 1    & 0.84 & 0.87 \\
 G3     & 1.15 & 0.93  & 0.97   & 1.67 & 0.85 & 1.18 & 1.39 & 1.07 & 1    & 0.93 \\
 G6     & 1.27 & 1.10  & 0.92   & 1.79 & 1.16 & 1.70 & 1.91 & 1.18 & 1.06 & 1    \\
  \hline
    \end{tabular*}

  \begin{tabular*}{\columnwidth}{@{}*{11}{c@{\extracolsep\fill}}c@{}}
    \hline
        & L    & EP1.5 & EP1.75 & N    & EP 3 & EP4  & U    & G1   & G3   & G6   \\
 L      & 2    & 1.12  & 1.02   & 1    & 1.07 & 1.15 & 1.34 & 1.46 & 1.53 & 1.18 \\
 EP1.5  & 1.16 & 2     & 1.02   & 1    & 0.87 & 0.55 & 0.40 & 1.03 & 1.26 & 1.89 \\
 EP1.75 & 1.03 & 1.21  & 2      & 1    & 0.90 & 0.83 & 0.81 & 1.00 & 1.04 & 1.15 \\
 N      & 1    & 1     & 1      & --   & 1    & 1    & 1    & 1    & 1    & 1    \\
 EP3    & 1.19 & 2.44  & 1.12   & 1    & 2    & 1.13 & 1.01 & 1.02 & 1.17 & 1.72 \\
 EP4    & 1.42 & 1.17  & 1.04   & 1    & 1.42 & 2    & 1.23 & 1.04 & 1.35 & 2.60 \\
 U      & 2.24 & 1.02  & 1.01   & 1    & 1.14 & 1.34 & 2    & 1.08 & 1.83 & 0.21 \\
 G1     & 2.32 & 1.18  & 1.03   & 1    & 1.19 & 1.74 & 2.24 & 2    & 1.20 & 1.05 \\
 G3     & 1.12 & 2.29  & 1.24   & 1    & 2.54 & 0.81 & 0.81 & 1.17 & 2    & 1.39 \\
 G6     & 1.04 & 1.15  & 1.91   & 1    & 0.94 & 0.93 & 0.94 & 1.05 & 1.29 & 2    \\

   \hline
    \end{tabular*}
\end{table}

\subsection{Finite-sample efficiencies}
\label{simulation:sec}

We compare the finite-sample efficiencies of the estimates in a simulation study using the same two-dimensional
settings with $\bo\Omega=\bo I_p$ as in the previous section. In each setting we consider the average of
$n(p-1)\hat{D}^2$ which has limiting expected value $ASV(\hat{\gamma}_{jl})+ASV(\hat{\gamma}_{lj})$. Thus, the
simulation study also illustrates how well the asymptotic results approximate the finite-sample variances. Let $\hat{\bo\Gamma}_i^{s2}$ and $\hat{\bo\Gamma}_i^{s}$, $i=1,\dots,M$,  be the estimates from $M$ samples of size $n$. Then the finite sample asymptotic relative efficiency is estimated by
\[
 \widehat{ARE}_{s2,s}  =\frac{\sum_{i=1}^M \{D(\hat{\bo\Gamma}_i^{s}\bo\Omega)^2\}}{\sum_{i=1}^M \{D(\hat{\bo\Gamma}_i^{s2}\bo\Omega)^2\}}.
\]

In Table~\ref{table2}, we list the estimated values of ${ARE}_{s2,s}$ and ${ARE}_{s2,d}$ for the same set of distributions as in
Table~\ref{table1}. For each setting, $M=10000$ samples of size $n=1000$ are generated. In most of the settings, the ratios
of the averages are close to the corresponding asymptotical values. When both components are nearly Gaussian, a larger
sample size than 1000 is required for $\widehat{ARE}_{s2,s}$     and $\widehat{ARE}_{s2,d}$  to converge to ${ARE}_{s2,s}$     and ${ARE}_{s2,d}$,
respectively. Also, if $\E[G(z_{ij})]\approx \E[G(z_{il})]$, then the extraction order of the deflation-based estimate
is not always the one which is assumed when computing the asymptotical variances. This may have a large impact on the
efficiency of the deflation-based estimate.

\begin{table*}[htb]
\renewcommand{\arraystretch}{1.3}
\caption{Values of $\widehat{ARE}_{s2,s}$  (on the top) and $\widehat{ARE}_{s2,d}$ (on the bottom)   computed from 10000 samples of size $n=1000$ for different distributions. L=Laplace=EP1=exponential power distribution with $\beta=1$, N=Normal distribution=EP2, U=Uniform distribution, G1=Gamma distribution with $\alpha=1$. Upper triangle for $pow3$ and lower triangle for $tanh$.}
\label{table2}
\centering
  \begin{tabular*}{16cm}{ccccccccccc}
    \hline
        & L    & EP1.5     & EP1.75             & N    & EP 3 & EP4  & U    & G1   & G3   & G6   \\
 L      & 0.98\textbackslash 0.86 & 0.94 & 1.09 & 1.17 & 0.80 & 0.75 & 0.73 & 0.85 & 0.88 & 0.98 \\
 EP1.5  & 1.03 & 1.01\textbackslash 0.95 & 1.15 & 1.34 & 0.94 & 1.07 & 1.49 & 0.92 & 0.89 & 0.95 \\
 EP1.75 & 1.27 & 1.27 & 1.14\textbackslash 1.13 & 1.07 & 1.66 & 2.10 & 3.13 & 1.02 & 1.06 & 1.08 \\
 N      & 1.47 & 1.65 & 1.12 & --                      & 2.22 & 2.99 & 4.13 & 1.05 & 1.25 & 1.25 \\
 EP3    & 0.92 & 1.06 & 1.59 & 1.91 & 1.31\textbackslash 1.07 & 1.08 & 1.45 & 0.85 & 0.77 & 0.90 \\
 EP4    & 0.84 & 1.14 & 1.67 & 2.12 & 1.06 & 1.07\textbackslash 1.00 & 1.14 & 0.81 & 0.76 & 0.98 \\
 U      & 0.87 & 1.41 & 1.89 & 2.22 & 1.25 & 1.08 & 1.00\textbackslash 1.00 & 0.76 & 0.82 & 1.23 \\
 G1     & 0.96 & 0.97 & 1.18 & 1.43 & 0.84 & 0.76 & 0.76 & 0.99\textbackslash 0.84 & 0.88 & 0.93 \\
 G3     & 1.15 & 0.94 & 1.12 & 1.51 & 0.87 & 1.01 & 1.38 & 1.13 & 1.10\textbackslash 0.83 & 0.90 \\
 G6     & 1.31 & 1.21 & 1.09 & 1.23 & 1.19 & 1.47 & 1.86 & 1.25 & 1.15 & 1.11\textbackslash 0.93  \\
  \hline
    \end{tabular*}

 \begin{tabular*}{16cm}{ccccccccccc}
   \hline
        & L    & EP1.5     & EP1.75             & N    & EP 3 & EP4  & U    & G1   & G3   & G6   \\
 L      & 1.93\textbackslash 1.79 & 1.12 & 1.02 & 1.00 & 1.09 & 1.18 & 1.38 & 1.47 & 1.53 & 1.18 \\
 EP1.5  & 1.14 & 1.82\textbackslash 1.47 & 1.14 & 0.98 & 1.66 & 1.54 & 0.85 & 1.04 & 1.29 & 1.40 \\
 EP1.75 & 1.03 & 1.24 & 1.14\textbackslash 1.04 & 1.03 & 1.18 & 0.89 & 0.72 & 1.00 & 1.04 & 1.11 \\
 N      & 1.00 & 0.98 & 1.05 & --                      & 0.89 & 0.91 & 0.92 & 1.00 & 1.00 & 1.00 \\
 EP3    & 1.18 & 1.78 & 1.19 & 0.93 & 1.82\textbackslash 1.78 & 1.35 & 1.01 & 1.03 & 1.25 & 1.50 \\
 EP4    & 1.42 & 1.41 & 1.01 & 0.97 & 1.40 & 1.93\textbackslash 1.94 & 1.23 & 1.05 & 1.42 & 1.52 \\
 U      & 2.14 & 1.00 & 0.99 & 0.99 & 1.13 & 1.33 & 1.92\textbackslash 1.94 & 1.11 & 1.64 & 1.30 \\
 G1     & 2.23 & 1.18 & 1.03 & 1.00 & 1.20 & 1.44 & 2.27 & 2.04\textbackslash 1.66 & 1.22 & 1.06 \\
 G3     & 1.43 & 1.99 & 1.20 & 0.99 & 1.79 & 1.85 & 1.03 & 1.26 & 2.18\textbackslash 1.57 & 1.36 \\
 G6     & 1.05 & 1.83 & 1.21 & 0.99 & 1.35 & 1.05 & 0.91 & 1.06 & 1.49 & 1.29\textbackslash 1.36 \\
   \hline
    \end{tabular*}
\end{table*}

\begin{figure*}[!t]
\centering
\includegraphics[width=5truecm]{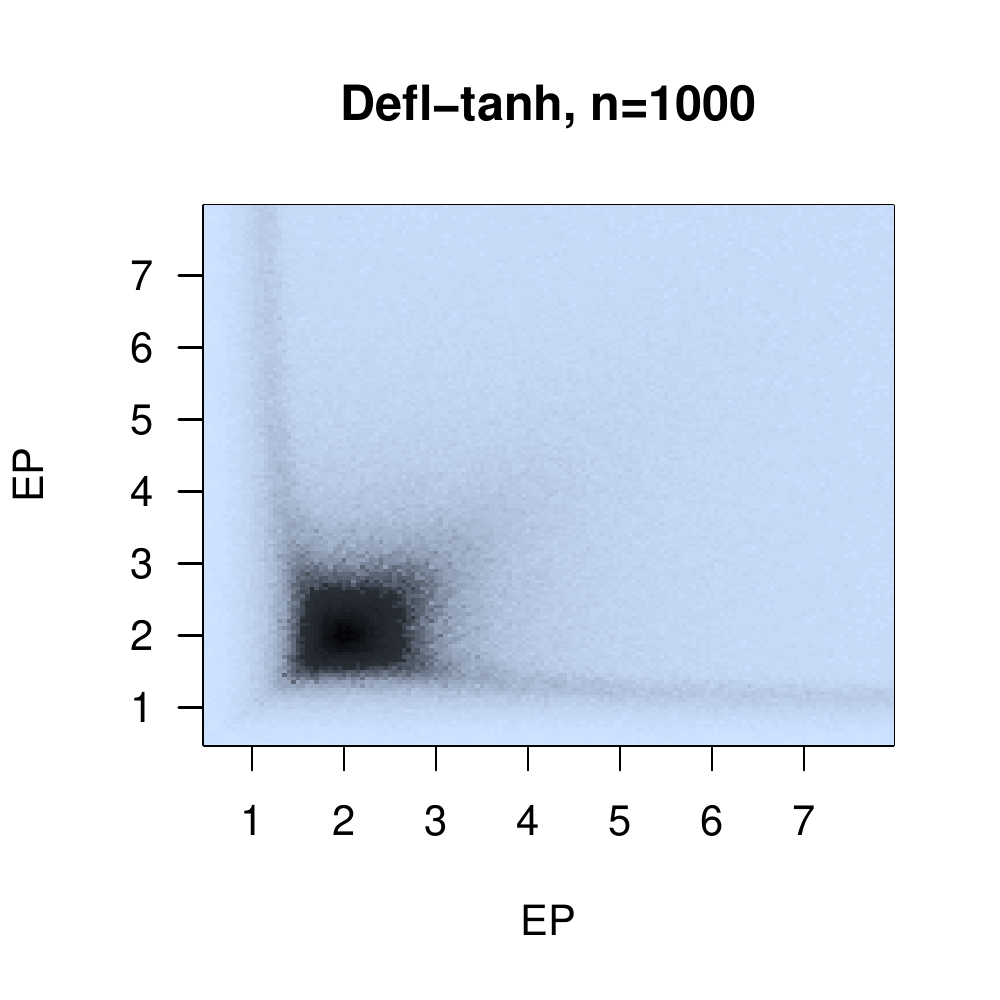}
\includegraphics[width=5truecm]{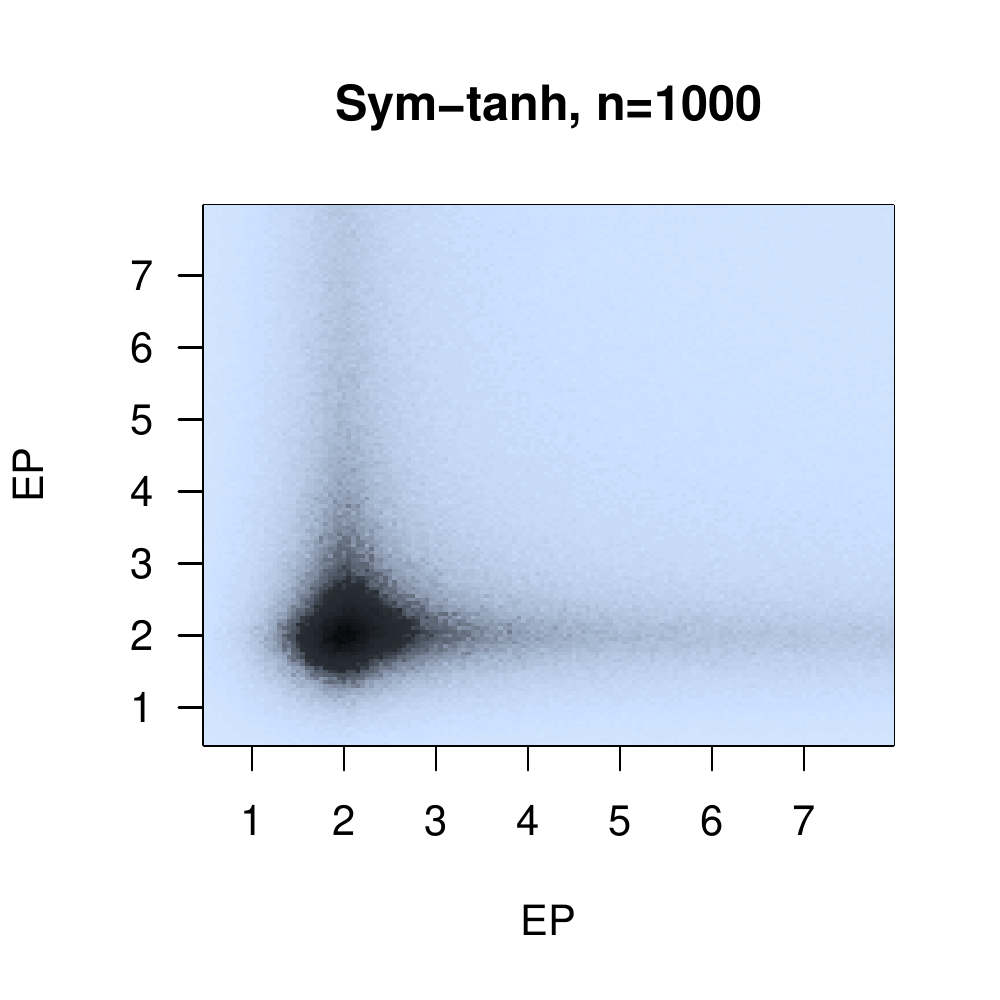}
\includegraphics[width=5truecm]{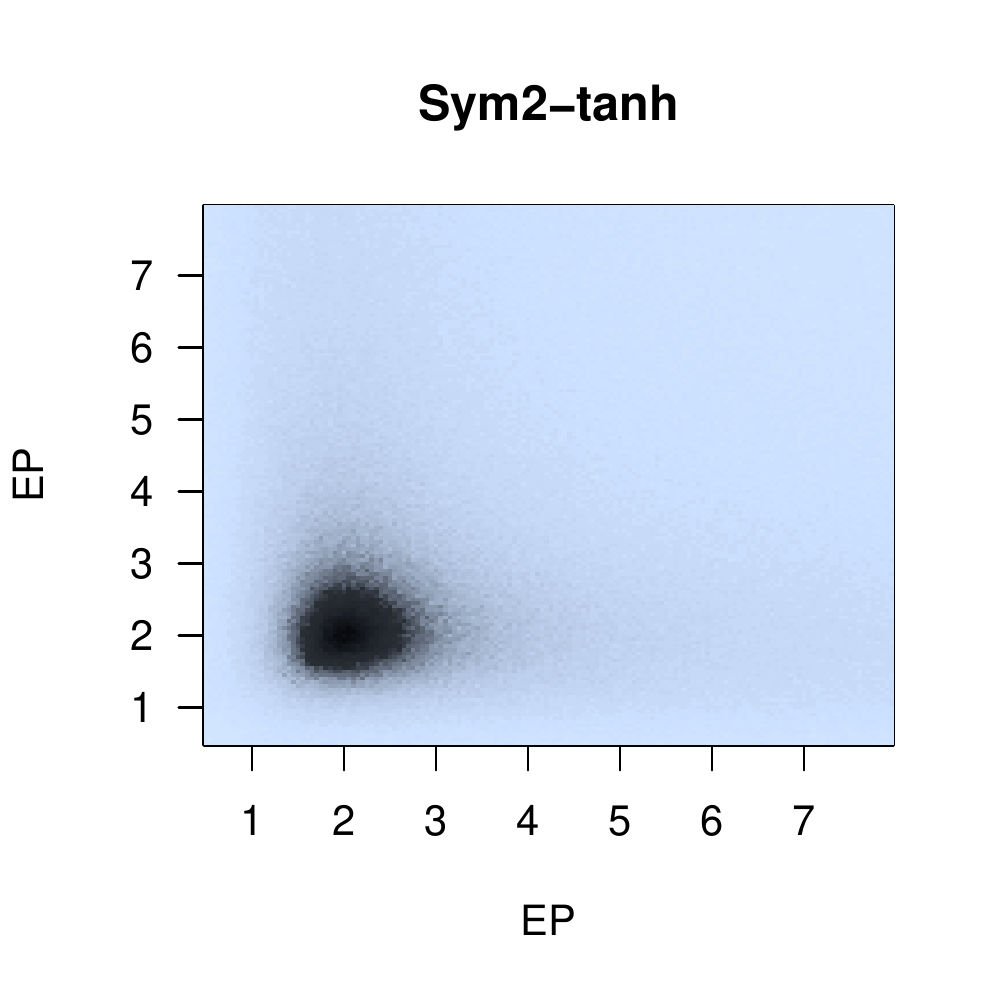}

\includegraphics[width=5truecm]{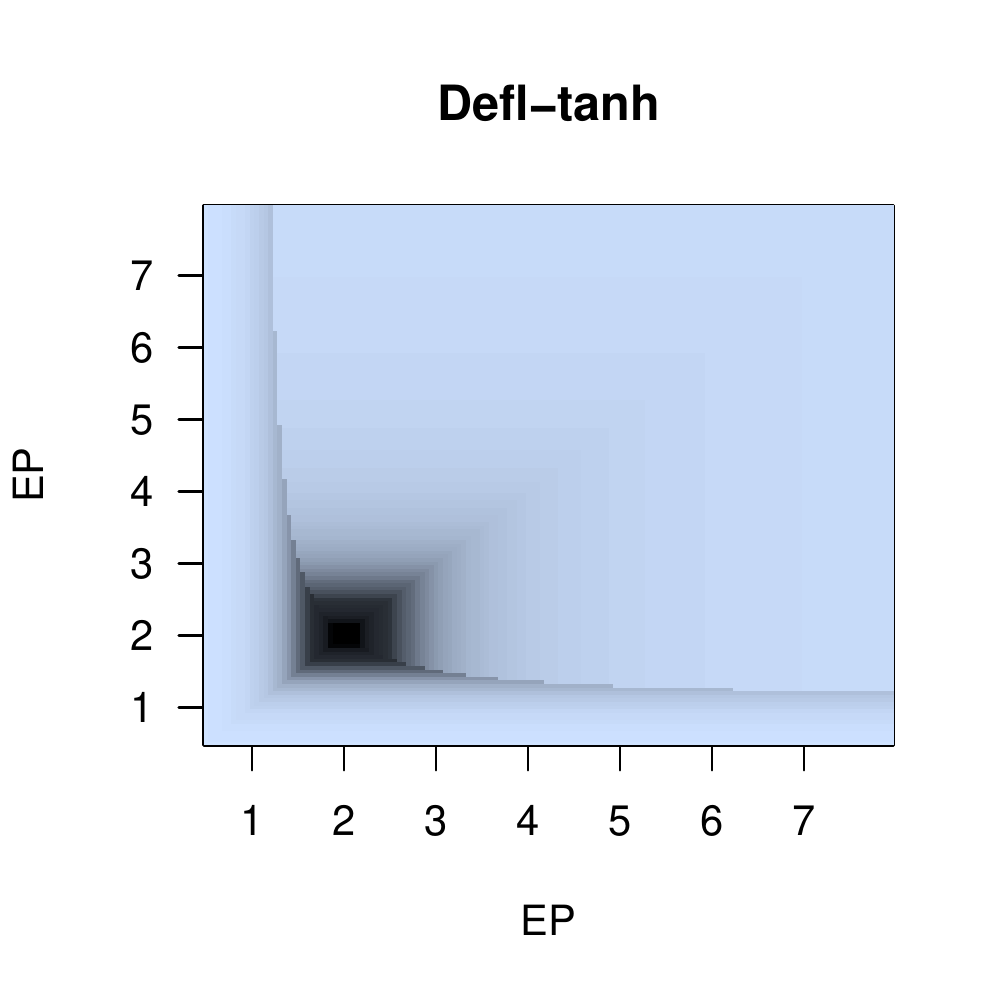}
\includegraphics[width=5truecm]{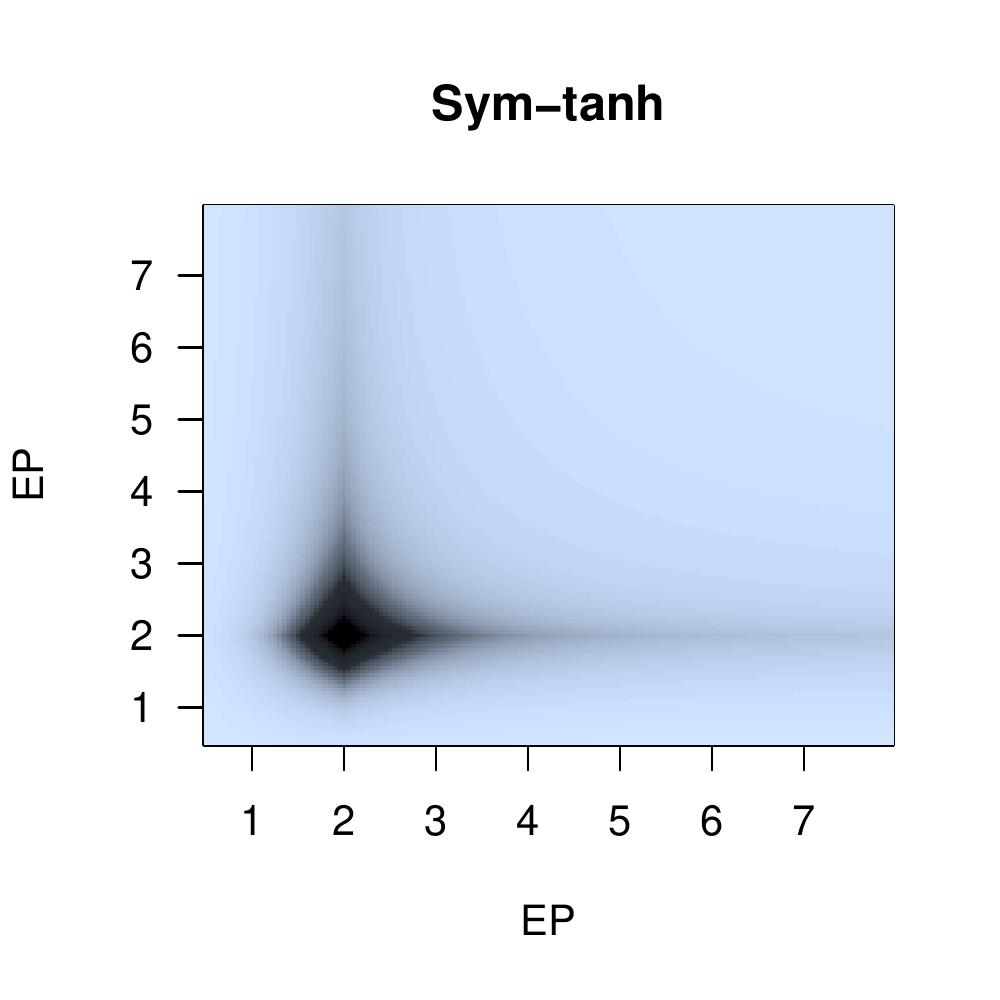}
\includegraphics[width=5truecm]{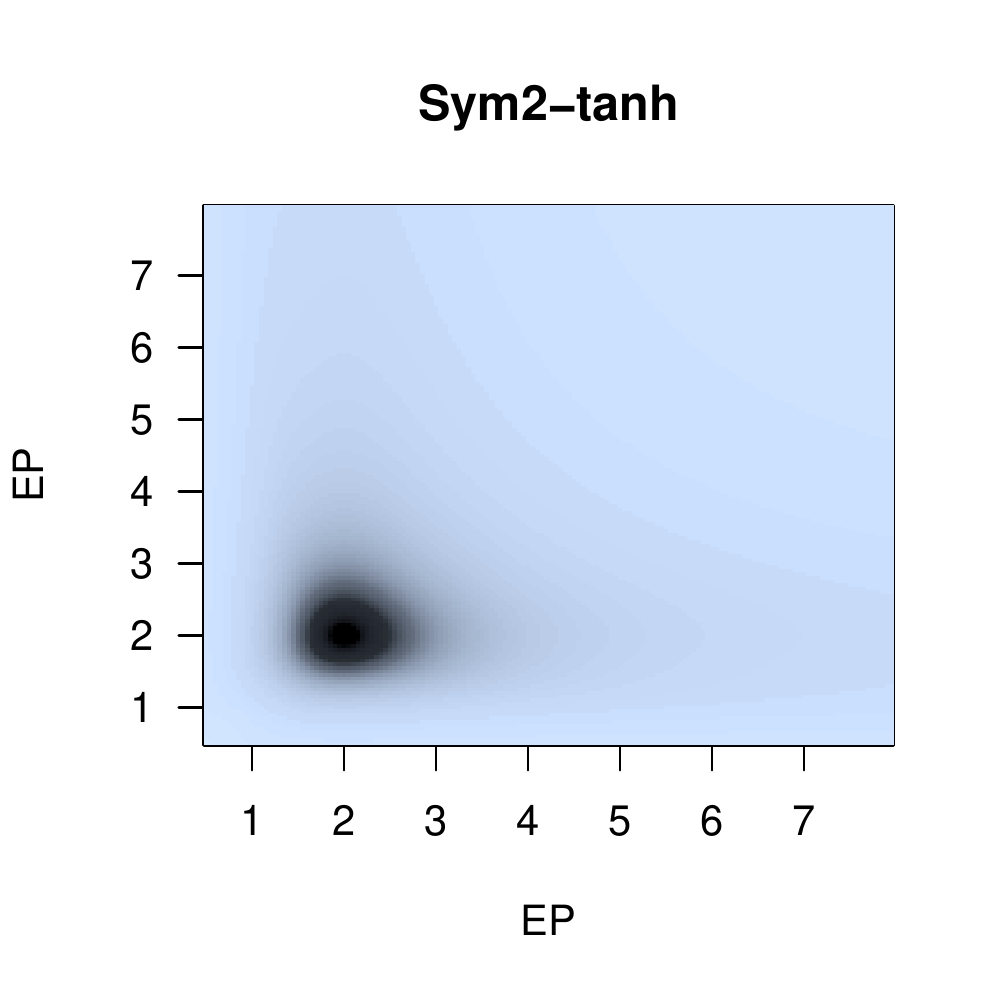}

\caption{Contour maps of the average of $n(p-1)\hat{D^2}$ over 200 simulation runs with
deflation-based, symmetric and squared symmetric FastICA estimates using {\it tanh} on the top and
the contour maps of the limiting expected values on the bottom. Two independent components with
exponential power distribution and varying shape parameter value. }
\label{fig4}
\end{figure*}

In Figure~\ref{fig4} we plot the contour maps of the average of $n(p-1)\hat{D}^2$ over 200 simulation runs
for deflation-based, symmetric and squared symmetric FastICA estimates using {\it tanh}. Each setting has two
independent components with exponential power distribution and varying shape parameter value, and $n=1000$.
Also the contour maps of the limiting expected values are given, and the corresponding maps resemble each other
rather nicely. The asymptotical results thus provide good approximations already for $n=1000$.

\section{Conclusions}

In this paper we investigate in detail the properties of the squared symmetric FastICA procedure, obtained from the regular symmetric FastICA procedure by replacing in the objective function the analytically cumbersome absolute values by their squares. We reviewed in a unified way the estimating equations, algorithms and asymptotic theory of the classical deflation-based and symmetric FastICA estimators and provided similar tools and derived similar results for the novel squared symmetric FastICA.  The asymptotic variances were used to compare the three methods in numerous different situations.

The asymptotic and finite sample efficiency studies imply, that although none of the methods uniformly outperforms  the others, the squared symmetric approach has the best overall performance under the considered combinations of source distributions and nonlinearities. Also, a crude ranking order of (\emph{deflation-based, symmetric, squared symmetric}) from worst to best can be given  and thus the use of the squared symmetric variant over the two other methods is highly recommended.

\appendices
\section*{Appendix}


Write
\begin{align*}
\hat{\bo T}(\hat{\bo\gamma})&=\frac{1}{n}\sum_{i=1}^n g(\hat{\bo\gamma}^T(\bo x_i-\bar{\bo x}))(\bo x_i-\bar{\bo x}) \text{ and} \\
\hat{\bo T}_2(\hat{\bo\gamma})&=\frac{1}{n}\sum_{i=1}^n G(\hat{\bo\gamma}^T(\bo x_i-\bar{\bo x}))\frac{1}{n}\sum_{i=1}^n g(\hat{\bo\gamma}^T(\bo x_i-\bar{\bo x}))(\bo x_i-\bar{\bo x}).
\end{align*}
The deflation-based, symmetric and squared symmetric FastICA estimators
$\hat{\bo\Gamma}^d=\bo\Gamma^d(\bo X)$, $\hat{\bo\Gamma}^s=\bo\Gamma^s(\bo X)$ and $\hat{\bo\Gamma}^{s2}=\bo\Gamma^{s2}(\bo X)$
are defined as follows
\begin{definition}
The deflation-based FastICA estimate $\hat{\bo\Gamma}^d=(\hat{\bo\gamma}_1^d,\dots,\hat{\bo\gamma}_p^d)^T$
solves the estimating equations
\begin{align}\label{defl2}
\hat{\bo T}(\hat{\bo\gamma}_j)=\hat{\bo S}\left(\sum_{l=1}^j\hat{\bo\gamma}_l\hat{\bo\gamma}_l^T\right)\,\hat{\bo T}(\hat{\bo\gamma}_j),\ \ \ j=1,\dots,p,
\end{align}
\end{definition}

\begin{definition}
The symmetric FastICA unmixing matrix estimate $\hat{\bo\Gamma}^s=(\hat{\bo\gamma}_1^s,\dots,\hat{\bo\gamma}_p^s)^T$
solves the estimating equations
\begin{align}\label{sym2}
\hat{\bo\gamma}_l^T\hat{\bo T}(\hat{\bo\gamma}_j)\,\hat s_j=\hat{\bo\gamma}_j^T\hat{\bo T}(\hat{\bo\gamma}_l)\,\hat s_l\ \ \ \text{and}\ \ \ \
\hat{\bo\gamma}^T_l\hat{\bo S}\hat{\bo\gamma}_j=\delta_{lj},
\end{align}
where $j,l=1,\dots,p$ and $\delta_{lj}$ is the Kronecker delta.
\end{definition}

\begin{definition}
The squared symmetric FastICA unmixing matrix estimate
$\hat{\bo\Gamma}^{s2}=(\hat{\bo\gamma}_1^{s2},\dots,\hat{\bo\gamma}_p^{s2})^T$ solves the estimating equations
\begin{align}\label{sqsym2}
\hat{\bo\gamma}_l^T\hat{\bo T}_2(\hat{\bo\gamma}_j)=\hat{\bo\gamma}_j^T\hat{\bo T}_2(\hat{\bo\gamma}_l)\ \ \ \text{and}\ \ \ \
\hat{\bo\gamma}^T_l\hat{\bo S}\hat{\bo\gamma}_j=\delta_{lj},
\end{align}
where $j,l=1,\dots,p$ and $\delta_{lj}$ is the Kronecker delta.
\end{definition}

To prove Theorem~\ref{asymp_theorem}, we need the following straightforward result:
\begin{lemma}
\label{gSg}
The second set of estimating equations $\hat{\bo\gamma}^T_j\hat{\bo S}\hat{\bo\gamma}_l=\delta_{lj}$, $j,l=1,\dots,p$
yields to
$$
\sqrt{n}\,(\hat{\gamma}_{jj}-1)=-\frac{1}{2}\sqrt{n}\,(\hat{\bo S}-\bo I_p)_{jj}  + o_P(1)
$$
and
\begin{align}
\label{gjl}
\sqrt{n}\,\hat{\gamma}_{jl}+\sqrt{n}\,\hat{\gamma}_{lj}=-\sqrt{n}\,\hat{\bo S}_{jl} + o_P(1).
\end{align}
\end{lemma}

\medskip\noindent{\bf Proof of Theorem~\ref{asymp_theorem} (iii)}

Let us now consider the first set of estimating equations. To shorten the notations, write
$\hat{\bo T}_2(\hat{\bo\gamma}_j)=\hat{\bo T}_{2j}$. Now
\[
\sqrt{n}\,\hat{\bo\gamma}_l^T\hat{\bo T}_{2j}=\sqrt{n}\,(\hat{\bo\gamma}_l-\bo e_l)^T\hat{\bo T}_{2j}+\sqrt{n}\,\bo e_l^T(\bo{\hat T}_{2j}-\nu_j\lambda_j\bo e_j).
\]
By Taylor expansion and Slutky's Theorem, we have
\begin{align*}
&\sqrt{n}\,(\hat{\bo T}_{2j}-\nu_j\lambda_j\bo e_j)=\sqrt{n}\,(\bo{T}_{2j}-\nu_j\lambda_j\bo e_j) \\
&-(\mu_j\lambda_j+\nu_j\tau_j)\bo e_j\bo e_j^T\sqrt{n}\,\bo{\bar x} \\ &+(\lambda_j^2\bo e_j\bo e_j^T+\nu_j\Delta_j)\sqrt{n}\,(\hat{\bo\gamma}_j-\bo e_j)+o_P(1),
\end{align*}
where $\Delta_j=\E[g'(z_{ij})\bo z_i\bo z_i^T]$. Consequently,
\begin{align*}
&\sqrt{n}\,\hat{\bo\gamma}_l^T\hat{\bo T}_{2j}=\sqrt{n}\,(\hat{\bo\gamma}_l-\bo e_l)^T\nu_j\lambda_j\bo e_j
 \\ & +\bo e_l^T(\sqrt{n}\,\bo{T}_{2j}-(\mu_j\lambda_j+\nu_j\tau_j)
\bo e_j\bo e_j^T\sqrt{n}\,\bo{\bar x} \\
& +(\lambda_j^2\bo e_j\bo e_j^T+\nu_j\Delta_j)
\sqrt{n}\,(\hat{\bo\gamma}_j-\bo e_j))+o_P(1) \\
&=\nu_j\lambda_j\sqrt{n}\,\hat{\gamma}_{lj}+\bo e_l^T\sqrt{n}\,\bo{T}_{2j}+\nu_j\delta_j\sqrt{n}\,\hat{\bo\gamma}_{jl}+o_P(1).
\end{align*}
According to our estimating equations, above expression should be equivalent to
\[
\sqrt{n}\,\hat{\bo\gamma}_j^T\hat{\bo T}_{2l}=\nu_l\lambda_l\sqrt{n}\,\hat{\gamma}_{jl}+\sqrt{n}\,\bo e_j^T\bo{T}_{2l}+\nu_l\delta_l\sqrt{n}\,\hat{\bo\gamma}_{lj}+o_P(1),
\]
which means that
\begin{align*}
&(\nu_l\lambda_l-\nu_j\delta_j)\sqrt{n}\,\hat{\gamma}_{jl}-(\nu_j\lambda_j-\nu_l\delta_l)\sqrt{n}\,\hat{\gamma}_{lj} \\
&=\sqrt{n}\,\bo e_l^T\bo T_{2j}-\sqrt{n}\,\bo e_j^T\bo T_{2l}+o_P(1).
\end{align*}
Now using~(\ref{gjl}) in Lemma~\ref{gSg}, we have that
\begin{align*}
&(\nu_l\lambda_l-\nu_j\delta_j)\sqrt{n}\,\hat{\gamma}_{jl}+(\nu_j\lambda_j-\nu_l\delta_l)\\
&(\sqrt{n}\,\hat{\bo S}_{jl}+\sqrt{n}\,\hat{\gamma}_{jl}) =\sqrt{n}\,(\bo e_l^T\bo T_{2j}-\bo e_j^T\bo T_{2l})+o_P(1).
\end{align*}
Then
\begin{align*}
&(\nu_j(\lambda_j-\delta_j)+\nu_l(\lambda_l-\delta_l))\sqrt{n}\,\hat{\gamma}_{jl} \\
=\ &\sqrt{n}\,(\bo e_l^T\bo T_{2j}-\bo e_j^T\bo T_{2l})+(\nu_l\delta_l-\nu_j\lambda_j)\sqrt{n}\,\hat{\bo S}_{jl}+o_P(1),
\end{align*}
which proves the Theorem.

The densities of $z_1$ and $z_2$ in Section~\ref{Gfunction:sec} are given by
\[
 f_i=\sum_{j=1}^4 \pi_{ij}N(\mu_{ij},\sigma_{ij}^2),\ \ i=1,2,
\]
where $N(\mu, \sigma^2)$ denotes the Gaussian density function with mean $\mu$ and variance $\sigma^2$,
and the (rounded) parameter values are
\begin{align*}
&\pi_{11}=0.09,\ \ \ \pi_{12}=0.43,\ \ \ \ \pi_{13}=0.43,\ \pi_{14}=0.04, \\
&\pi_{21}=0.15,\ \ \ \pi_{22}=0.31,\ \ \ \ \pi_{23}=0.45,\ \pi_{24}=0.09, \\
&\mu_{11}=-1.76,\ \mu_{12}=-0.34,\ \mu_{13}=0.54,\ \mu_{14}=1.79, \\
&\mu_{21}=-1.71,\ \mu_{22}=-0.36,\ \mu_{23}=0.48,\ \mu_{24}=1.66, \\
&\sigma_{11}^2=0.13,\ \ \ \sigma_{12}^2=0.50,\ \ \ \ \sigma_{13}^2=0.28,\ \sigma_{14}^2=0.13, \\
&\sigma_{21}^2=0.11,\ \ \ \sigma_{22}^2=0.26,\ \ \ \ \sigma_{23}^2=0.11,\ \sigma_{24}^2=0.11.
\end{align*}

\section*{Acknowledgment}

This work was supported by the Academy of Finland (grants 251965, 256291 and 268703).

\ifCLASSOPTIONcaptionsoff
  \newpage
\fi

\end{document}